\documentclass[12pt]{article}

\usepackage{amsmath, epsfig, cite, lineno}
\usepackage{amssymb,amsthm}
\usepackage{amsfonts, color}
\usepackage{latexsym}

\newtheorem{thm}{Theorem}[section]
\newtheorem{prop}[thm]{Proposition}

\newtheorem{conj}[thm]{Conjecture}
\newtheorem{lem}[thm]{Lemma}

\newcommand{\pf}{\noindent{\it Proof.} }

\numberwithin{equation}{section}

\begin{document}


\begin{center}
{\Large\bf Proof of two conjectures of Z.-W. Sun on\\[5pt] congruences for Franel numbers}
\end{center}

\vskip 2mm \centerline{Victor J. W. Guo}
\begin{center}
{\footnotesize Department of Mathematics, East China Normal
University,\\ Shanghai 200062,
 People's Republic of China\\
{\tt jwguo@math.ecnu.edu.cn,\quad
http://math.ecnu.edu.cn/\textasciitilde{jwguo} }}
\end{center}


\vskip 0.7cm \noindent{\bf Abstract.}
For all nonnegative integers $n$, the Franel numbers are defined as
\begin{align*}
f_n=\sum_{k=0}^n {n\choose k}^3.
\end{align*}
We confirm two  conjectures of Z.-W. Sun on congruences for
Franel numbers:
\begin{align*}
\sum_{k=0}^{n-1}(3k+2)(-1)^k f_k  &\equiv 0 \pmod{2n^2},\\
\sum_{k=0}^{p-1}(3k+2)(-1)^k f_k  &\equiv 2p^2 (2^p-1)^2 \pmod{p^5},
\end{align*}
where $n$ is a positive integer and $p>3$ is a prime.



\section{Introduction}

The numbers $f_n$ are defined to be the sums of cubes of binomial coefficients:
\begin{align*}
f_n=\sum_{k=0}^n {n\choose k}^3.
\end{align*}
In 1894, Franel \cite{Franel1,Franel2} obtained the following recurrence relation for $f_n$:
\begin{align}
(n+1)^2 f_{n+1}=(7n^2+7n+2)f_{n}+8n^2 f_{n-1},\ n=1,2,\ldots. \label{eq:JV}
\end{align}
Nowadays, the numbers $f_n$ are usually called Franel numbers. The Franel numbers also appear
in the first and second Strehl identities \cite{Strehl1,Strehl2} (see also Koepf \cite[p.~55]{Koepf}):
\begin{align*}
\sum_{k=0}^n {n\choose k}^3 &=\sum_{k=0}^n{n\choose k}^2{2k\choose n}, \\
\sum_{k=0}^n {n\choose k}^2 {n+k\choose k}^2
&=\sum_{k=0}^n\sum_{j=0}^k {n\choose k}{n+k\choose k}{k\choose j}^3.
\end{align*}

Applying the recurrence relation \eqref{eq:JV}, Jarvis  and Verrill \cite{JV} proved the following congruence for Franel numbers
$$
f_{n}\equiv (-8)^n f_{p-1-n} \pmod p,
$$
where $p$ is a prime and $0\leq n\leq p-1$.
Recently, Z.-W. Sun \cite{Sun}, among other things, proved several interesting congruences for
Franel numbers, such as
\begin{align*}
\sum_{k=0}^{p-1}(-1)^{k} f_k &\equiv \left(\frac{p}{3}\right) \pmod{p^2}, \\
\sum_{k=0}^{p-1}(-1)^{k} k f_k &\equiv -\frac{2}{3}\left(\frac{p}{3}\right) \pmod{p^2},\\
\sum_{k=0}^{p-1}(-1)^{k} k^2 f_k &\equiv \frac{10}{27}\left(\frac{p}{3}\right) \pmod{p^2},
\end{align*}
where $p>3$ is a prime and $\left(\frac{a}{3}\right)$ denotes the Legendre symbol. Sun \cite{Sun} also
proposed many amazing conjectures on congruences for $f_n$. The main purpose of this  paper is
to prove the following results, which were conjectured by Sun \cite{Sun}.

\begin{thm}\label{thm:2nn} For any positive integer $n$, there holds
\begin{align}
\sum_{k=0}^{n-1}(3k+2)(-1)^k f_k \equiv 0 \pmod{2n^2}. \label{eq:f3}
\end{align}
\end{thm}

\begin{thm}\label{thm:ppower5} For any prime $p>3$, there holds
\begin{align}
\sum_{k=0}^{p-1}(3k+2)(-1)^k f_k \equiv 2p^2 (2^p-1)^2 \pmod{p^5}.
\end{align}
\end{thm}

\section{Proof of Theorem \ref{thm:2nn}}
We need the following identity due to MacMahon \cite[p.~122]{MacMahon}
(see also Foata \cite{Foata} or Riordan \cite[p.~41]{Riordan}):
\begin{align}
\sum_{k=0}^{n}{n\choose k}^3x^k
=\sum_{k=0}^{\lfloor n/2\rfloor}{n+k\choose 3k}{3k\choose 2k}{2k\choose k}x^k (1+x)^{n-2k}.
\label{eq:foata}
\end{align}
When $x=1$, the above identity \eqref{eq:foata} gives a new expression for Franel numbers:
\begin{align}
\sum_{k=0}^{n}{n\choose k}^3
=\sum_{k=0}^{\lfloor n/2\rfloor}{n+k\choose 3k}{3k\choose 2k}{2k\choose k}2^{n-2k}. \label{eq:mac}
\end{align}
Differentiating both sides of \eqref{eq:foata} twice with respect to $x$ and then letting $x=1$, we get
\begin{align}
\sum_{k=0}^{n}{n\choose k}^3 k(k-1)
=\sum_{k=0}^{\lfloor n/2\rfloor}2^{n-2k}{n+k\choose 3k}{3k\choose 2k}{2k\choose k}\frac{n(n-1)-2k}{4}.
\label{eq:mac2}
\end{align}
Moreover, by induction, we can easily prove that
\begin{align}
\sum_{\ell=2k}^{n-1}(-1)^\ell (3\ell +2){\ell+k\choose 3k}2^{\ell-2k}
=(-1)^{n-1}(n-2k){n+k\choose 3k}2^{n-2k}. \label{eq:induc}
\end{align}

In fact, when $n=2k$, both sides of \eqref{eq:induc} are equal to $0$. Now suppose that \eqref{eq:induc}
is true for $n$. Then
\begin{align*}
&\hskip -2mm \sum_{\ell=2k}^{n}(-1)^\ell (3\ell +2){\ell+k\choose 3k}2^{\ell-2k} \\
&=(-1)^n (3n+2){n+k\choose 3k}2^{n-2k}+\sum_{\ell=2k}^{n-1}(-1)^\ell (3\ell +2){\ell+k\choose 3k}2^{\ell-2k}\\
&=(-1)^n (3n+2){n+k\choose 3k}2^{n-2k}+(-1)^{n-1}(n-2k){n+k\choose 3k}2^{n-2k} \\
&=(-1)^{n}(n-2k+1){n+k+1\choose 3k}2^{n-2k+1},
\end{align*}
which implies that \eqref{eq:induc} holds for $n+1$.

Applying \eqref{eq:mac} and then exchanging the summation order, we have
\begin{align}
\sum_{k=0}^{n-1}(3k+2)(-1)^k f_k
&=(-1)^{n-1}\sum_{k=0}^{n-1}2^{n-2k}(n-2k){n+k\choose 3k}{3k\choose 2k}{2k\choose k} \label{eq:sum1}
\end{align}
in view of \eqref{eq:induc}. Noticing that
$$
n-2k =4 \frac{n(n-1)-2k}{4} -(n^2-2n),
$$
by \eqref{eq:mac} and \eqref{eq:mac2}, we can write the right-hand side of \eqref{eq:sum1} as
\begin{align*}
&\hskip -3mm (-1)^{n-1}4\sum_{k=0}^{n}{n\choose k}^3 k(k-1)+(-1)^{n}(n^2-2n)\sum_{k=0}^{n}{n\choose k}^3\\
&=(-1)^{n-1}4n^2\sum_{k=0}^{n}{n\choose k}{n-1\choose k-1}^2 +(-1)^n n^2\sum_{k=0}^{n}{n\choose k}^3,
\end{align*}
where we have used the following relations:
\begin{align*}
k{n\choose k}&=n{n-1\choose k-1}, \\
\sum_{k=0}^{n}k{n\choose k}^3&=\sum_{k=0}^{n}(n-k){n\choose k}^3=\frac{n}{2}\sum_{k=0}^{n}{n\choose k}^3.
\end{align*}
Namely, we have proved that
\begin{align}
\frac{1}{2n^2}\sum_{k=0}^{n-1}(3k+2)(-1)^k f_k
&=(-1)^{n-1}2\sum_{k=0}^{n}{n\choose k}{n-1\choose k-1}^2 +(-1)^n \frac{f_n}{2}.
\label{eq:sum2}
\end{align}
The proof then follows from the fact
$$
f_n=\sum_{k=0}^{n}{n\choose k}^3\equiv \sum_{k=0}^{n}{n\choose k}=2^n \equiv 0 \pmod 2,\  n\geq 1.
$$

By \eqref{eq:mac}, for $n\geq 1$, we have
\begin{align*}
f_n\equiv
\begin{cases}
2,&\text{if $n$ is a power of $2$,}\\[5pt]
0,&\text{otherwise,}
\end{cases}  \pmod 4.
\end{align*}
In fact, if $n=2m+1\geq 3$ is odd, then ${2k\choose k}2^{2m-2k+1}\equiv 0\pmod 4$ for all
$k\leq m$ and so $f_{2m+1}\equiv 0\pmod 4$; if $n=2m$ is even,
then $f_{2m}\equiv{3m\choose 2m}{2m\choose m}=2{3m\choose 2m}{2m-1\choose m}\pmod 4$
and the result follows from the congruences:
\begin{align*}
{2m-1\choose m}\equiv
\begin{cases}
1,&\text{if $m$ is a power of $2$,}\\[5pt]
0,&\text{otherwise,}
\end{cases}  \pmod 2
\end{align*}
and ${3m\choose 2m}\equiv 1\pmod 2$ if $m$ is a power of $2$.

Thus, we may further refine Theorem \ref{thm:2nn} as follows:
\begin{thm}\label{thm:2nnnew} For any positive integer $n$, there holds
\begin{align}
\sum_{k=0}^{n-1}(3k+2)(-1)^k f_k
\equiv
\begin{cases}
\displaystyle 2n^2,&\text{if $n$ is a power of $2$,}\\[5pt]
0,&\text{otherwise,}
\end{cases} \pmod{4n^2}.
\end{align}
\end{thm}

\section{Proof of Theorem \ref{thm:ppower5}}

The following lemma is due to Sun \cite[Lemma 2.1]{Sun}.
\begin{lem}[Sun]\label{lem:sun}
For any prime $p>3$, there holds
\begin{align}
f_{p-1}\equiv 1+3(2^{p-1}-1)+3(2^{p-1}-1)^2 \pmod{p^3}. \label{eq:sun}
\end{align}
\end{lem}

To prove Theorem \ref{thm:ppower5}, we also need the following variation of Lemma \ref{lem:sun}.

\begin{lem}
For any prime $p>3$, there holds
\begin{align}
\sum_{k=1}^{p-1}{p-1\choose k}{p-1\choose k-1}^2
\equiv 2^{p-1}-2^{2p-2} \pmod{p^3}. \label{eq:sunvar}
\end{align}
\end{lem}
\pf  It easily follows from cubing ${p\choose k}={p-1\choose k}+{p-1\choose k-1}$ that
\begin{align}
\sum_{k=0}^{p}{p-1\choose k}^3
&=\sum_{k=0}^{p}{p\choose k}^3 -\sum_{k=0}^{p}{p-1\choose k-1}^3
-3\sum_{k=0}^{p}{p\choose k}{p-1\choose k}{p-1\choose k-1}. \label{eq:fpsum}
\end{align}
 Since
${p\choose k}\equiv 0\pmod p$ for $0<k<p$, we have
\begin{align}
\sum_{k=0}^{p}{p\choose k}^3\equiv 2\pmod{p^3}. \label{eq:fpmod}
\end{align}
Substituting \eqref{eq:sun} and \eqref{eq:fpmod} into \eqref{eq:fpsum},
we immediately get
\begin{align}
\sum_{k=0}^{p}{p\choose k}{p-1\choose k}{p-1\choose k-1}
\equiv 2^p-2^{2p-1}\pmod{p^3}. \label{eq:sun2}
\end{align}

On the other hand, replacing $k$ by $p-k$, we obtain
\begin{align}
\sum_{k=1}^{p-1}{p-1\choose k}{p-1\choose k-1}^2
= \sum_{k=1}^{p-1}{p-1\choose k}^2{p-1\choose k-1}
=\frac{1}{2}\sum_{k=0}^{p}{p\choose k}{p-1\choose k}{p-1\choose k-1}. \label{eq:sun3}
\end{align}
Combining \eqref{eq:sun2} and \eqref{eq:sun3}, complete the proof.  \qed

\medskip
\noindent{\it Proof of Theorem \ref{thm:ppower5}.}
By \eqref{eq:sum2} and \eqref{eq:fpmod}, we have
\begin{align*}
\frac{1}{2p^2}\sum_{k=0}^{p-1}(3k+2)(-1)^k f_k
&\equiv 2\sum_{k=0}^{n}{p\choose k}{p-1\choose k-1}^2 - 1,  \\
&=2\sum_{k=0}^{n}{p-1\choose k}{p-1\choose k-1}^2+2\sum_{k=0}^{n}{p-1\choose k-1}^3 - 1.
\end{align*}
The proof then follows from \eqref{eq:sun} and \eqref{eq:sunvar}.  \qed

\section{Concluding remarks and open problems}
For any nonnegative integers $n$ and $r$, let
$$
f_{n}^{(r)}=\sum_{k=0}^n{n\choose k}^r.
$$
Then $f_{n}^{(3)}=f_n$ are the Franel numbers. Calkin \cite[Proposition 3]{Calkin} proved the following congruence:
\begin{align*}
f_{n}^{(2r)}  \equiv 0 \pmod{p},
\end{align*}
where $p$ is a prime such that $\frac{n}{m}<p<\frac{n+1}{m}+\frac{n+1-m}{m(2mr-1)}$ for some positive integer $m$.
Guo and Zeng \cite[Theorem 4.4]{GZ} proved  that,  for any positive integer $n$,
\begin{align}
f_{n}^{(2r)}\equiv 0 \pmod {n+1}. \label{eq:gz}
\end{align}
Sun \cite[Conjecture 3.5]{Sunconj} conjectured that
\begin{align}
\sum_{k=0}^{n-1}(3k+2)f_{k}^{(4)}\equiv 0\pmod {2n}. \label{eq:f4}
\end{align}
It is easy to see that $f_{n}^{(0)}=n+1$, $f_{n}^{(1)}=2^n$, $f_{n}^{(2)}={2n\choose n}$
by the Chu-Vandermonde identity (see \cite[p.~41]{Koepf}), and
\begin{align}
\sum_{k=0}^{n-1}(3k+2)f_{k}^{(0)} &=n^3+n^2, \label{eq:f0}\\
\sum_{k=0}^{n-1}(-1)^k (3k+2)f_{k}^{(1)} &=(-1)^{n-1}2^n n,\label{eq:f1}\\
\sum_{k=0}^{n-1}(3k+2)f_{k}^{(2)} &=n{2n\choose n}. \label{eq:f2}
\end{align}
Motivated by the identities \eqref{eq:f0}--\eqref{eq:f2}, the congruence \eqref{eq:f3}, and
Sun's conjecture \eqref{eq:f4}, we would like to propose the following conjecture on congruences for $f_n^{(r)}$.
\begin{conj} \label{conj1}
Let $n\geq 1$ and $r\geq 0$ be two integers. Then
\begin{align}
\sum_{k=0}^{n-1}(-1)^{rk}(3k+2)f_{k}^{(r)}\equiv 0\pmod {2n}. \label{eq:fn}
\end{align}
\end{conj}
By \eqref{eq:gz}, if the congruence \eqref{eq:fn} holds, then we have
\begin{align*}
\sum_{k=0}^{n-1}(3k+2)f_{k}^{(2r)}\equiv 0\pmod {n(n+1)}. 
\end{align*}

For example, the first values of $\sum_{k=0}^{n-1}(3k+2)f_{k}^{(6)}$ are
$$
2, 12, 540, 16600, 784500, 35315784, 1772807064, 90283679280, 4777960538340,
$$
while $\sum_{k=0}^{n-1}(-1)^{k}(3k+2)f_{k}^{(5)}$ gives
$$
2,-8,264,-5104,132460,-3373824,91312256,-2513335808,70719559668.
$$

It seems that, for $n>1$, the following congruence holds:
\begin{align*}
\sum_{k=0}^{n-1}(-1)^{k}(3k+2)f_{k}^{(2r+1)}\equiv 0\pmod {4n}.
\end{align*}

Recall that the multinomial coefficients are given by
$$
{n\choose k_1,\ldots,k_m}=\frac{n!}{k_1!\cdots k_m!},
$$
where $k_1,\ldots,k_m\geq 0$ and $k_1+\cdots+k_m=n$.
 Let
$$
M_{m,n}^{(r)}
=\sum_{k_1+\cdots+k_m=n}{n\choose k_1,\ldots,k_m}^r
$$
be the sums of $r$th powers of multinomial coefficients.
Then $M_{m,n}^{(0)}={m+n-1\choose n}$, $M_{m,n}^{(1)}=m^r$, $M_{1,n}^{(r)}=1$, $M_{2,n}^{(r)}=f_{n}^{(r)}$,
and
$$
M_{3,n}^{(2)}=\sum_{k=0}^n \sum_{j=0}^{k} {n\choose k}^2 {k\choose j}^2=\sum_{k=0}^{n}{n\choose k}^2{2k\choose k}.
$$
Note that Osburn and Sahu \cite{OS} studied supercongruences for the numbers $\sum_{k=0}^{n}{n\choose k}^r{2k\choose k}^s.$
The sequence $\{M_{3,n}^{(2)}\}_{n\geq 0}$ is
the A002893 sequence of Sloane \cite{Sloane}. It also appears in Zagier \cite[\#8 of Table 1]{Zagier}.
Sun\cite{Sun} proved the following identity involving $M_{3,k}^{(2)}$:
\begin{align}
\sum_{k=0}^{n-1}(4k+3)M_{3,k}^{(2)}=3n^2\sum_{k=0}^{n-1}\frac{1}{k+1}{2k\choose k}{n-1\choose k}^2. \label{eq:sunmk3}
\end{align}
It seems that Conjecture \ref{conj1} can be further generalized as follows.
\begin{conj}\label{conj2}
Let $m,n\geq 1$ and $r\geq 0$ be integers. Then
\begin{align}
\sum_{k=0}^{n-1}(-1)^{rk}((m+1)k+m)M_{m,k}^{(r)}\equiv 0\pmod {mn}. \label{eq:mnr}
\end{align}
\end{conj}

It is not hard to prove the following identities by induction.
\begin{align}
\sum_{k=0}^{n-1}((m+1)k+m)M_{m,k}^{(0)} &=mn{m+n-1\choose m}, \label{eq:mnr1}\\
\sum_{k=0}^{n-1}(-1)^{k}((m+1)k+m)M_{m,k}^{(1)} &=(-1)^{n-1}m^n n, \label{eq:mnr2}\\
\sum_{k=0}^{n-1}(2k+1)M_{1,k}^{(2r)} &=n^2 , \label{eq:mnr3} \\
\sum_{k=0}^{n-1}(-1)^{k}(2k+1)M_{1,k}^{(2r+1)} &=(-1)^{n-1}n. \label{eq:mnr4}
\end{align}
Combining the congruence \eqref{eq:f3}, the identities \eqref{eq:f2},
\eqref{eq:sunmk3} and \eqref{eq:mnr1}--\eqref{eq:mnr4}, we have the following result to
support Conjecture \ref{conj2}.

\begin{prop}The congruence \eqref{eq:mnr} holds if $(m,r)$ belongs to
$$
\{(m,0),(m,1),(1,r),(2,2),(2,3),(3,2)\colon m,r\geq 1\}.
$$
\end{prop}


\vskip 5mm
\noindent{\bf Acknowledgments.}  This work was partially
supported by the Fundamental Research Funds for the Central Universities.

\end{document}